\documentclass[11pt]{article}
\usepackage{amssymb}
\usepackage{latexsym}
\usepackage{amsthm}

\newtheorem{de}{Definition}[section]

\newtheorem{te}[de]{Theorem}
\newtheorem{co}[de]{Corollary}
\newtheorem{exs}[de]{Examples}

\newtheorem{os}[de]{Remark}
\newtheorem{lem}[de]{Lemma}

\newcommand{\bfn}{\mbox{{\sf I}\hspace{-0.4em}{\sf N}}}

\newcommand{\gen}[1]{\mbox{$\langle #1 \rangle$}}

\title{Coverings of the Symmetric and Alternating Groups\footnote{Work supported by
M.U.R.S.T. and G.N.S.A.G.A} }
\author{Daniela Bubboloni}
\date{}
\begin {document}
\maketitle
\begin{abstract}
We analyse for which $n$ there exist in $G=A_n,S_n$ two proper subgroups $H,K$
such that $G$ is the  union of the $G$-conjugacy classes of $H$ and $K$.
\end{abstract}

\section{Introduction}

Let $G$ be a finite group. A set $\delta$ of proper subgroups $H_1,\ldots ,H_n$ of $G$ is called
a {\em covering} of $G$ if $G$ is the set-theoretical union of the $H_i$ and moreover, to avoid
 redundance, there are no inclusions between the $H_i$:
\[G=\bigcup_{i=1}^{n}H_i,\quad\quad  H_i\not \leq H_j\quad \hbox{for } i\neq j.\]

 We refer to the $H_i$ as the {\em components} of the covering. From the well known fact that
 {\em a group is never the union of two proper subgroups} it follows that {\em a covering  of a group  needs
 at least three components}.
A group  $G$ for which  a cover exists is called {\em coverable}.
It is clear that each group  admitting a partition is coverable, hence we have soon a lot of
available examples.
On the other hand, here, there are no assumption about the intersection of the components, so
 that the great part of the arguments used in the partition's theory becomes useless.
This is also one of the reason for the difficulty to develop a general and, at the same time,
 expressive description of the coverable groups. Not by accident, this theme appears in the
recent literature always from a particular point of view.
Brandl {\cite{bra81}} considers the coverings $\delta=\{H^{\alpha}: \alpha\in Aut(G)\}$ where
$H\leq G$; Praeger in {\cite{pra88} and {\cite{pra94}} explores the more general coverings
$\delta=\{H^{\alpha}:\alpha\in A\}$ where $Inn(G)\leq A\leq Aut(G)$ to study the coverings of
 the Galois group of certain extensions of fields.

We believe that a reasonable idea in exploring the coverable groups is to begin with some
 "natural" sets of subgroups which cannot be used to cover a group, trying to add to them some
"natural" new components in order to obtain a covering. In this sense our starting point is a
 well known fact: {\em  a finite group  $G$ is never the set-theoretical union of the
$G$-conjugates of a proper subgroup $H$}. Adding to the
$\,G$-conjugates of $H$ another subgroup $K$ or the set of the
$\,G$-conjugates of a subgroup $K,$ we have the opportunity, in
some cases, to obtain a covering; this happens, for example, in a
Frobenius group if $H$ and $K$ are respectively a complement and
the kernel. In other words we are led to consider two kinds of
collections of subgroups for a group $G$:
\[(\ast)\;\delta=\{H^g,K:g\in G\},\]

\[(\ast\ast)\;\delta=\{H^g,K^g:g\in G\}\]
where $H$ and $K$ are fixed proper subgroups of $G$. Observe that if $\delta$ is of the type
$(\ast)$ or $(\ast\ast)$ and $G$ is the set-theoretical union of the subgroups in $\delta,$ then
no inclusions among the elements of $\delta$ are possible, that is $\delta$ is a covering of $G.$
If a group $G$ is coverable with a set of subgroups of the type $(\ast)$ or $(\ast\ast)$, we
will say respectively that $G$ is
{\em $(\ast)$-coverable} or {\em $(\ast\ast)$-coverable}. In the next section we give a
characterization of the $(\ast)$-coverable groups, from which we deduce, in particular,
that no simple group is $(\ast)$-coverable. But it is easily seen that there exist simple
$(\ast\ast)$-coverable groups, such as $A_5$; the  whole family of the alternating groups
looks  interesting on this regard and we devote to this subject the section $3$, showing that
$A_n$ is $(\ast\ast)$-coverable if and only if  $\;4\leq n\leq 8$. In the last section
we obtain a  similar result for the symmetric group proving that $S_n$ is
$(\ast\ast)$-coverable if and only if  $\;3\leq n\leq 6$.

All groups in this paper are finite. All unexplained notation is standard ( see
\cite{hup67} ).

This paper has greatly benefited from very helpful discussions with C. Casolo.\\
I would also like to thank S. Dolfi regard to \ref{3.4} and P. Neumann regard to the proof
of \ref{2.10}. A special thank to E. De Tomi for the help in preparing this final version.

\bigskip

\section{$(\ast)$-coverable groups }

In this section, we aim to describe the $(\ast)$-coverable groups $G,$ that is the groups
$G$ coverable with a set of subgroups
of the type \[(\ast)\;\delta=\{H^g,K:g\in G\},\] where $H$ and $K$ are fixed proper subgroups
of $G$.

We begin with some elementary facts.

\begin{lem}\label{2.1}

Let  $G$ be a  group and $\delta=\{H^g,K:g\in G\}$ a covering of $G$.Then:
\begin {enumerate}
\item[(i)] $H^G=G$ and $H_G \leq K$;
\item[(ii)]if $K\lhd G$, then $G=KH$.
\end{enumerate}
\end{lem}
\begin{proof}
(i) Since $G=K \cup \left( \bigcup_{g\in G}H^g \right)$, it follows that $G=K\cup H^G$ with
$K\neq G$ hence $G=H^G$.
Assume now to have $x \in H_G$  with $x\notin K$. If we pick any $\,y\in K$, then we get
$xy\notin K$ hence there exists $g\in G$ such that $xy\in H^g$. But $x\in H^g$ and
therefore $y\in H^g$. This means that $K \subseteq \bigcup_{g\in G}H^g$, hence
$G=\bigcup_{g\in G}H^g$ and $H=G$, a contradiction.
\\
(ii)  If $K\lhd G$, then we simply observe that $G=\bigcup_{g\in G}KH^g=\bigcup_{g\in G}(KH)^g$
concluding $G=KH$.
\end{proof}

The next result shows that there is no loss of generality in assuming $K$ normal in $G$.

\begin{lem} \label{2.3}

Let  $G$ be a  group and $\delta=\{H^g,K:g\in G\}$ a covering of $G$. Then also
$\hat{\delta}=\{K_G,H^g:g\in G\}$ is a covering of $G$.
\end{lem}
\begin{proof}
Let $g\in G-\bigcup_{a\in G}H^a$. Then $\delta$ covering yields $g\in K$ and if $y\in G$, then
we necessarily have $g^{y^{-1}}\in K$, otherwise there would exist an element $x\in G$ with
$g^{y^{-1}}\in H^x$, that is $g\in H^{xy}$, against the assumption on $g$. Therefore for
every $y\in G$ we have $g^{y^{-1}}\in K$, that is $g\in K^y$ and then $g\in K_G$.\\
Clearly $K_G\neq G$ because $K\neq G$ by definition of covering and hence
$\hat{\delta}$ is a covering.
\end{proof}
\begin{co}\label{2.4}
No simple group is $(\ast)$-coverable.
\end{co}

To exhibit some crucial examples of $(\ast)$-coverable groups, we need to recall the definition
 of the Frobenius-Wielandt groups and the fundamental theorem related to them.

\begin{de}{\em (Wielandt){\cite{wi58}}}\label{2.5}
{\em A group $G$ is said to be a} Frobenius-Wielandt group {\em provided that it has a subgroup $H$,
 with $1\neq H\neq G$, and a proper normal subgroup $N$ of $H$ such that}
\end{de}
\[H \cap H^g\leq N\quad \mbox{if}\;g\in G-H.\]

For brevity, we will refer to them as {\em F-W groups} and we will use the locution
 {\em $(G,H,N)$ is a F-W group} to indicate more closely the situation.

\begin{te}{\em (Wielandt)\cite {wi58}}\label{2.6}
If $(G,H,N)$ is a F-W group, then there exists a unique normal subgroup $K$ of $G$ ( called the
 $\underline{kernel}$ ) such that:
\[G=HK, \quad H\cap K=N.\]
Moreover $K=G-\bigcup_{g\in G}(H-N)^g$.
\end{te}
\begin{lem}\label{2.7}
Let $(G,H,N)$ be a F-W group  with kernel $K$. Then \[\delta=\{H^g,K:g\in G\}\] is a
covering of $G$.
\end{lem}
\begin{proof}
In fact if $x\in G-K$, by \ref{2.6},  we get $x\in \bigcup_{g\in G}(H-N)^g$ hence $x\in H^g$
for some $g\in G$. Moreover \ref{2.6} and \ref{2.5} yield that $H,K\neq G,$ hence
$\delta=\{H^g,K:g\in G\}$ is a covering of $G$.
\end{proof}

It is clear that each Frobenius group $G$ with complement $H$ is a $(G,H,1)$  F-W group.
Furthermore, as a trivial consequence of the definition \ref{2.5}, we get that each group with
a F-W quotient is itself a F-W group. Particularly, we get:
\begin{os}\label{2.8}
If a group has a Frobenius quotient, then it is a F-W group.
\end{os}

This fact allows us to be more concrete in the construction of examples of $(\ast)$-coverable
groups.
\begin{exs}\label{2.9}

$S_3$, $A_4$ and $S_4$ are $(\ast)$-coverable groups.
\end{exs}
\begin{proof}
In fact, by \ref{2.8}, these are all F-W groups and then \ref{2.7} applies.
\end{proof}
Now an application of the counting method leads to the converse of \ref{2.7}.
\begin{lem}\label{2.10}
If the group $G$ has a covering $\delta=\{H^g,K:g\in G\}$, then $(G,H,H\cap K_G)$ is a F-W
group with kernel $K_G$.
\end{lem}
\begin{proof}
Due to \ref{2.3}, we may assume that $K\lhd G$. First of all we show that $H$ is
selfnormalizing in $G$. Let $H_1,\ldots ,H_n$ be the distinct conjugates of $H$ in $G$ and
$m=|N_G(H):H|$. Then we have $|G|=mn|H|<\sum_{i=1}^n|H_i| +|K|$, because
$G=\left(\bigcup_{i=1}^{n} H_i\right)\cup K$.
From $K<G$ we get also $\sum_{i=1}^n|H_i| +|K|\leq n|H|+\frac{|G|}{2}=n|H|+
\frac{mn}{2}|H|$, hence $mn<n\left(1+\frac{m}{2} \right)$, that is $m<2$ and
then $m=1$.

Next we define the subgroups $K_i=H_i\cap K$ and we prove that the subsets $X_i=H_i-K_i$  have
empty intersection for $i\neq j$. By \ref{2.1}(ii), we have $G=KH_i$ hence
$|G|=\frac{|K||H|}{|K_i|}$ that is $|K_i|=\frac{|K|}{n}$ for each $i=1,\ldots,n$. It
follows that $|K|+\sum_{i=1}^n|X_i|=\sum_{i=1}^n|K_i|+\sum_{i=1}^n\left(|H_i|-|K_i|
\right)=n|H|=|G|$.
On the other hand, obviously, we have $G=K\cup \left( \bigcup_{i=1}^{n} X_i\right)$ and then
the previous relation implies $X_i\cap X_j=\emptyset$ for $i\neq j$. Clearly this means that $H_i\cap
H_j \leq K$. Let $g\in G-H=G-N_G(H)$, then $H^g\neq H$ and therefore $H^g\cap H\leq H\cap K$,
that is $(G,H,H\cap K)$ is a F-W group. Finally from the characterization of the kernel of
a F-W group given in \ref{2.6}, it follows that $K$ is actually the  kernel of $G$.
\end{proof}
Collecting \ref{2.7} and \ref{2.10}, we can state the main result of this section
\begin{te}\label{2.11}
A group $G$ is $(\ast)$-coverable if and only if $G$ is a F-W group.
\end{te}
\bigskip

\section{The $(\ast\ast)$-coverable alternating groups}

Given a set $\Omega$, we denote the symmetric and the alternating groups on
$\Omega$ respectively with {\em $Sym\,\Omega$} and {\em $Alt\,\Omega$}.
When $\Omega=\{1,\ldots,n\}$ we use, more simply, the notations {\em $S_n$} and {\em $A_n$}
and we consider the natural immersions of $S_n$ into $S_{n+1}$ and of $A_n$ into $A_{n+1}$ as
inclusions.

If $\sigma \in S_n$ decomposes into the product of disjoint cycles $\sigma_1,\ldots,\sigma_k$
of lengths $l_1,\ldots,l_k$ we will say that the {\em type} of $\sigma$ is {\em $[l_1;
\ldots;l_k]$}. Sometimes, when not misleading, we will omit the lengths equal to $1$.\\
Obviously it holds:
\begin{os}\label{3.1}
Let $G=S_n$ and $H,\;K<G$.Then $\delta=\{H^g,K^g:g\in G\}$ is a covering of $G$ if and only if
each type of permutation appears at least one time either in $H$ or in $K$.
\end{os}

Though no simple group is $(\ast)$-coverable ( \ref{2.4} ), it looks reasonable to investigate the
simple groups $G$ which admit a $(\ast\ast)$-covering that is which are covered by a set of
subgroups of the type $(\ast\ast)\;\delta=\{H^g,K^g:g\in G\}$ for some $H,\ K<G.$
There is in fact a natural example: $A_5$ is
covered by \{$A_4^g, P^g : g\in A_5$\} where $P\in Syl_5(G)$, just because
$\bigcup_{g\in A_5}A_4^g$ contains all the permutations in $A_5$ with at least a fixed
point and $\bigcup_{g\in A_5}P^g$ all the permutations with no fixed point. On the other
hand this construction is peculiar for $A_5$ and not extendible to the alternating groups
of higher degree. Thus it is not evident for which $n\in \bfn$, $A_n$ is
$(\ast\ast)$-coverable.\\
Observe that if $H,\ K<A_n$ and each type of even permutation lies in $H$ or in $K$,
this does not guarantee any more that $\delta=\{H^g,K^g:g\in A_n\}$ is a covering of $A_n$
( cf. \ref{3.1} ). The problem is of course that certain types of even permutations
decompose into two conjugacy classes in $A_n.$ We need to be more specific about this
question:
\begin{lem}\label{3.4}
If $\sigma \in A_n$, then the $S_n$-conjugacy class $\sigma^{S_n}$ splits into two
$A_n$-conjugacy classes if and only if $\sigma$ is of the type $[l_1;\ldots;l_r]$ with
$l_i\geq 1$ distinct in pairs and odd for $i=1,\ldots,r.$
\end{lem}
\begin{proof} {\em ( probably folklore )} First of all we observe that if $\sigma \in A_n,$ then
$\sigma^{S_n}$ splits into two $A_n$-conjugacy classes if and only if $C_{S_n}(\sigma)\leq A_n.$

Next we note that if $\mu \in S_n$ is a cycle, then
\begin{equation}
C_{S_n}(\mu)= \gen{\mu}\times Sym(\Omega - supp(\mu)) 
\end{equation}
where, for each $\sigma \in S_n,$ we denote with {\em supp($\sigma$)} ( the support of $\sigma$ )
the set of $i\in \Omega=\{ 1,\ldots,n\}$ such that $i^{\sigma}\neq i.$ Namely it is enough to
show $(1)$ for the $k$-cycle $\mu=(1\,2\ldots k)$ with $2\leq k\leq n$ and obviously
$\gen{\mu}\times Sym(\Omega - supp(\mu))\leq C_{S_n}(\mu).$
Moreover if $\alpha \in C_{S_n}(\mu),$ then we can write $\alpha=
\overline{\alpha}\beta$ where $\overline{\alpha}$ is a product of disjoint cycles whose support
has non empty intersection with $K=\{1,\ldots,k\}$ and $\beta$ is a product of disjoint cycles
whose support has empty intersection with $K.$ Since $\alpha$ stabilizes $K,$ we get that $supp(
\overline{\alpha}) \subseteq K,$ $\beta \in Sym\{k+1,\ldots,n\}$ and $\overline{\alpha}\in
C_{S_n}(\mu).$ Therefore we have $(1\,2\ldots k)^{\overline{\alpha}}=(1^{\overline{\alpha}}\,2^
{\overline{\alpha}}\ldots k^{\overline{\alpha}})=(1\,2\ldots k)$ and then, if we put
$1^{\overline{\alpha}}=i\in K,$ we get $\,j^{\overline{\alpha}}=i+j-1\,(\,\hbox{mod}\ k\,)\,$
for each $1\leq j\leq k.$ But we have also $j^{\mu^{i-1}}=i+j-1\,(\,\hbox{mod}\ k\,)\,$ for each
$1\leq j\leq k,$ that is $\,\overline{\alpha}=\mu^{i-1} \in \gen{\mu}\,$ and so
$\alpha \in \gen{\mu}\times Sym(\Omega - supp(\mu)).$

Now we observe that if $\sigma \in S_n$ and $\sigma=\mu_1\cdots \mu_r$ with $\mu_i$ disjoint
cycles of lengths $l_1,\ldots,l_r$ distinct in pairs, then
\begin{equation}
C_{S_n}(\sigma)= \bigcap_{i=1}^{r}C_{S_n}(\mu_i).  
\end{equation}
In fact the inclusion $\bigcap_{i=1}^{r}C_{S_n}(\mu_i)\leq C_{S_n}(\sigma)$ is trivial;
moreover if $\alpha \in C_{S_n}(\sigma),$ then $\,\mu_1^{\alpha}\ldots \mu_r^{\alpha}=
\mu_1\ldots \mu_r\,$  are two decompositions in the product of disjoint cycles and
$\mu_i^{\alpha}$ is
a $\,l_i$-cycle. Since $\mu_i$ is the only $l_i$-cycle among the $\mu_1,\ldots ,\mu_r,$ we
argue that $\mu_i^{\alpha}=\mu_i$ for each $i=1,\ldots,r$ that is
$\alpha \in \bigcap_{i=1}^{r}C_{S_n}(\mu_i).$

Let $\sigma \in A_n$ be the product of the disjoint cycles $\mu_1,\ldots ,\mu_r,$ of lengths
$l_i\geq 1.$\\ If the $l_i$ are odd and distinct in pairs, then by $(2)$ and $(1)$ we obtain
\[C_{S_n}(\sigma)= \bigcap_{i=1}^{r}C_{S_n}(\mu_i)=\bigcap_{i=1}^{r}\left[\gen{\mu_i}\times
Sym(\Omega - supp(\mu_i))\right]=\times_{i=1}^{r}\gen{\mu_i}\leq A_n.\]
Therefore $\sigma^{S_n}$ splits into two $A_n$-conjugacy classes.\\ Now assume that there exist
at least two $\mu_i$ of the same odd length, say $l_1=l_2=l$ odd and let $\mu_1=(i_1\ldots i_l), \
\mu_2=(i_1'\ldots i_l').$ Then we have $(i_1\,i_1')\ldots (i_l\,i_l')\in C_{S_n}(\sigma)-A_n.$
Finally if at least one among the $\mu_i$ has even length, then clearly $\mu_i\in
C_{S_n}(\sigma)-A_n.$ In both cases we have $C_{S_n}(\sigma)\not \leq A_n$ and hence
$\sigma^{S_n}=\sigma^{A_n}.$
\end{proof}
We approach the problem of the determination of the $n\in \bfn$ for which $A_n$ is
$(\ast\ast)$-coverable, beginning with $n\leq 8$. To do this, first of all, we state an
elementary lemma which will be useful also in the sequel and analyse $S_5$ and $S_6.$
\begin{lem}\label{3.3}
Let $G$ be a group covered by $\delta=\{H^g,K^g:g\in G\}$. If $N\unlhd G$ and $G=NH=NK$,
then $N$ admits the covering $\delta_N=\{(H \cap N)^x,(K\cap N)^x\-:x\in N\}$.
\end{lem}

We refer to the covering {\em $\delta_N$} defined in the previous lemma as a {\em $(\ast\ast)$-covering
obtained by intersection}.
\begin{os}\label{3.2}
$S_5,\;S_6$ are $(\ast\ast)$-coverable subgroups.
\end{os}
\begin{proof}
Let $G=S_5,$\[ H_1=G_{\{1,2\}}\] and \[K_1=N_G\gen{(12345)}.\]Then $K_1=\gen{(12345)}\rtimes
\gen{(2354)}<G$
contains permutations of the types $[4],\ [5],\ [2;2]$ and $H_1<G$ contains permutations of
the types $[2],\ [3],\ [2;3]$. Hence, by \ref{3.1}, $\delta=\{H_1^g,K_1^g:g\in G\}$ is a
covering of $G$.

Now let $G=S_6$, \[H_2\hbox { the stabilizer in G of the partition }\{1,2,3\},\ \{4,5,6\}\]
and\[K_2=S_5.\]Then $H_2\cong  S_3\ wr\ S_2$ contains the permutation $(142536)$ of type $[6]$
and therefore also a permutation of the type $[2;2;2]$ and one of the type $[3;3]$. Moreover
$H_2$ contains the permutation $(14)(2536)$ of the type $[2;4]$; hence for each type of
fixed-point-free permutation $H_2$ contains at least a representative. On the other hand
$K_2$ contains representatives for each type of permutation with at least a fixed point and
again \ref{3.1} applies.
\end{proof}
\begin{os}\label{3.5}
$A_5,\;A_6,\;A_7$  and $A_8$ are $(\ast\ast)$-coverable groups.
\end{os}
\begin{proof}
We easily obtain a covering of $A_5$ and $A_6$ by intersection from the coverings constructed in
\ref{3.2}  for $S_5$ and $S_6$. In fact $A_n \lhd S_n$ and, in
order to apply \ref{3.3}, we only need to observe that  each of the subgroups $H_1,H_2,K_1,K_2$
contains an odd permutation. But it is trivially checked that
 \[(12)\in H_1,\ (2354)\in K_1,\  (12)\in H_2\hbox{ and }(23)\in K_2.\]

Next let $G=A_7$ and consider \[H=N_G\gen{(1234567)}\] and \[K=[Sym\{1,2\}\times
Sym\{3,4,5,6,7\}]\cap A_7.\] Since
$N_{S_7}\gen{(1234567)}=\gen{(1234567)}\rtimes \gen{\mu}$ with $\mu$ a $6$-cycle, we
get that $H=\gen{(1234567)}\rtimes\gen{\mu^2}<G$ where $\mu^2$ is of type $[1;3;3].$
By \ref{3.4}, the
permutations of type $[1;3;3]$ constitute a single conjugacy class hence
$\bigcup_{g\in G}H^g$ contains all the permutations of this type. Moreover $H$ contains
a $7$-Sylow of $G$ hence $ \bigcup_{g\in G}H^g$ contains every $7$-cycle of $G$.
On the other hand, from  $|K|_2=|G|_2$, it follows that
in $\bigcup_{g\in G}K^g$ there is each
$2$-element of $G$. Next observe that $K$ contains at least a permutation of type
$[3],\ [5],\ [2;2;3]$ and therefore, again by \ref{3.4}, each permutation of these
types lies in $\bigcup_{g\in G}K^g$. Since we have examined all the possible types of
permutations in $G$, then it follows that $\{H^g,K^g:g\in G\}$ is a covering of $G=A_7$.

Finally we explore the group $A_8.$ Let $H$ be the group of all the affine transformations
of the vector space $V=GF(2)^3,$ that is of all the transformations $\tau_{_{A,a}}
:V\rightarrow V$ of the form
\[\tau_{_{A,a}}(x)=Ax+a\]
where $A\in GL(3,2)$ and $a\in V.$\\
Obviously, $H\leq Sym\,V\cong S_8.$ We want to show that actually $H<Alt\,V\cong A_8.$ It is
well known that $H=T\rtimes GL(3,2),$ where $T$ is the elementary abelian subgroup of order
$2^3$ consisting of all the translations $\tau_{_{I,a}}\in H.$ Moreover, it is clear that each
$1\neq \tau \in T$ has no fixed points on $V,$ hence $\tau \in Sym\,V$ is a permutation of
the type $[2;2;2;2]$ and $T\leq Alt\,V.$ On the other hand, $GL(3,2)$ is simple and
therefore we necessarily have $GL(3,2) \leq Alt\,V,$ otherwise $GL(3,2)$ would contain a
normal subgroup of index two. Thus we have $H<Alt\,V$ and, because $|H|=2^6\cdot 3\cdot 7,$
it is clear that $H$ contains a $7$-cycle.\\
Next we observe that, putting
\[a=\left(\begin{array}{c}0\\1\\0\end{array}\right)\in V \hbox{ and }
A=\left(\begin{array}{ccc}1&1&0\\0&1&0\\0&0&1\end{array}\right)\in GL(3,2),\]
it is easily checked that $\tau_{_{A,a}}\in H$
is a permutation on $V$ of the type $[4;4].$\\
Moreover, putting
\[b=\left(\begin{array}{c}1\\0\\0\end{array}\right)\in V \hbox{ and }
B=\left(\begin{array}{ccc}1&0&0\\0&0&1\\0&1&1\end{array}\right)\in GL(3,2),\]
we get that $\tau_{_{B,b}}\in H$
is a permutation on $V$ of the type $[2;6].$\\
It follows that $H$ contains permutations of the types:
\[(1)\quad\quad [2;2;2;2],\ [4;4],\ [2;6],\ [7].\]
Therefore, by \ref{3.4} and by the Sylow theorem, we get that  $\bigcup_{g\in A_8}H^g$
contains every permutation of these types.\\
Now we denote with the natural numbers $1,\ldots,8$ the elements in $V,$ we identify
$Alt\,V$ with $A_8$ and we put
\[K=\left[Sym\,\{1,2,3\}\times Sym\,\{4,5,6,7,8\}\right]\cap A_8.\]
Thus $K$ contains at least a permutation of the types:
\[[3],\ [5],\ [2;2],\ [2;4],\ [3;3],\ [2;2;3]\]
and also the two non-conjugated permutations $(123)(45678)$ and $(123)^{(12)}(45678)$ of
the type $[3;5].$ Again, by \ref{3.4}, this implies that $\bigcup_{g\in A_8}K^g$ contains
all the permutations of the types
\[(2)\quad \quad[3],\ [5],\ [2;2],\ [2;4],\ [3;3],\ [3;5],\ [2;2;3].\]
Since each type of permutation in $A_8$ belongs either to the list $(1)$ or to the
list $(2),$ we conclude that $\{H^g,K^g:g\in A_8\}$ is a covering of $A_8.$
\end{proof}

Note that if $A_n$ is $(\ast\ast)$-coverable, {\em there is no sort of uniqueness for the
$(\ast\ast)$-coverings} of $A_n$: we have already showed two different
$(\ast\ast)$-coverings for $A_5$ and it is easily seen that we can construct another
$(\ast\ast)$-covering for $A_6$ using
\[H=Alt\{2,3,4,5,6\},\]\[K=\gen{(14)(2356),(15)(24)}\cong S_4.\]

At this point, to continue our investigation on the $(\ast\ast)$-coverable $A_n$, the
leading concept becomes the primitivity. In fact we will appeal
several times to the  following classical  result.
\begin{te}{\em ( \cite{wipe}, 13.9)}\label{3.6}  
Let $G\leq S_n$ be a primitive group and $p$ a prime such that $p\leq n-3.$ If $G$ contains
a $p$-cycle, then $G\geq A_n.$
\end{te}

We recall that if $G\leq S_n$ is $2$-transitive then $G$ is primitive and also that
if $G\leq S_p,$ with $p$ a prime, is transitive then $G$ is primitive.
Moreover, for our purpose, it will be often appropriate to deduce the primitivity of a
permutation group from the next lemma.
\begin{lem}\label{3.7}
Let $G\leq S_n$ be a transitive subgroup and $n_0$ the minimal non trivial divisor of $n.$ If there
exists a prime $p$ such that $p>\frac{n}{n_0}$ and $p\mid\,|G|$, then $G$ is primitive.
\end{lem}
\begin{proof}
Assume that $\Delta$ is a non trivial block of imprimitivity for $G$. By the transitivity of $G$,
$|\Delta|\mid \,n\,$ and therefore $n_0\leq |\Delta|\leq \frac {n}{n_0}.$

Let $\overline{\Omega}=\{\Delta_1=\Delta,\ldots,\Delta_l\}$ be a complete system of blocks
for $G$ and consider the action of $\gen{\sigma}$ on $\overline{\Omega}$, where $\sigma
\in G$ is an element of order $p.$
This action cannot be faithful, otherwise $Sym\,\overline{\Omega}$ would contain an element
of order $p$ that is the product of some disjoint $p$-cycles; consequently
$\,p\leq |\overline{\Omega}|=\frac{n}{|\Delta|}\leq \frac {n}{n_0}$ and we would reach a
contradiction.

This means that the action of $\gen{\sigma}$ on $\overline{\Omega}$ is trivial and if
$(i_1,\ldots,i_p)$ is one of the disjoint $p$-cycles in which $\sigma$ splits then,
renumbering the $\Delta_i$, we can assume that $i_1\in \Delta=\Delta^{\sigma}.$
In particular $\Delta\supseteq\{i_1,\ldots,i_p\}$ and then
$|\Delta|\geq p>\frac {n}{n_0}$ which gives again a contradiction.
\end{proof}
\begin{os} \label{trans}
Let $\{H^g,K^g:g\in A_n\}$ be a covering of $A_n.$ If $n\geq 5,$ then at least one among $H$
and $K$ is transitive. If $n\geq 7,$ then exactly one among $H$ and $K$ is transitive.
\end{os}
\begin{proof}
Let $\{H^g,K^g:g\in A_n\}$ be a covering of $A_n$ and $n\geq 5.$\\
We can assume that $H$ contains a cycle $\sigma$ of maximal length in $A_n.$ If
$n$ is odd, then $H$ is clearly transitive.\\
If $n$ is even $\sigma$ has length $n-1$ and if
$H$ has more than one orbit, then it has exactly two orbits of lengths $1$ and $n-1$.
In this case in $\bigcup_{g\in A_n}H^g$ we get exclusively permutations with at least a fixed
point and therefore $K$ contains a permutation $\lambda$ of the type $[2;n-2]$ and a
permutation $\mu$ of the type $[3;n-3]$. To fix the ideas  assume that $\lambda$
interchanges $1$ and $2$. Since $n\geq 6$, it follows that in the decomposition of $\mu$
as a product of disjoint cycles there are no transpositions hence $\mu$ takes $1$ or
$2$ into an element belonging to $\{3,\ldots ,n\}$. But $\lambda$ cyclically permutes
the elements of $\{3,\ldots ,n\}$ and consequently $K$ is transitive on
$\{1,\ldots ,n\}.$

Next let $n\geq 7$ and assume that $H$ and $K$ are both transitive. If $n=7,$ then $H$ and
$K$ are primitive and one of them
must contain a $3$-cycle hence, by \ref{3.6}, it coincides with $A_n$. Then we can
assume $n\geq 8$. By the Bertrand's postulate, there exists
a prime $p$ with $n/2<p\leq n-3$ and either $H$ or $K$ must contain a $p$-cycle
$\sigma$, say $\sigma \in H.$ By \ref{3.7}, we get that $H$ is primitive and, from
\ref{3.6}, it follows that $H=A_n,$ a contradiction.
\end{proof}

We are now in position to prove the goal of this section.
\begin{te}\label{3.9}
$A_n$ is  $(\ast\ast)$-coverable if and only if $\,4\leq n\leq 8.$
\end{te}
\begin{proof}
From \ref{2.9} and \ref{3.5} we know that if $4\leq n\leq 8,$ then $A_n$ is a
$(\ast\ast)$-coverable group.
Thanks to \ref{trans} we only need to show that for $n>8,$ $A_n$ admits no covering
$\{H^g,K^g:g\in A_n\}$ with $H,\ K<A_n,\ H$ transitive and $K$ not transitive.
Assume the contrary. First of all we observe that
if $[l_1;\ldots;l_r]$ is the type of a permutation in $A_n$, then at least a
permutation of this type lies in $H$ or in $K.$\\
Let $\gamma_1,\ldots,\gamma_k$, with $k\geq 2$ be the orbits of $K$ on $\Omega=
\{1,\ldots,n\}$ and $\gamma_1$ that of minimal length $a$. Then
\[K\leq\left[Sym\,\gamma_1 \times Sym\,(\gamma_2\cup\ldots\cup\gamma_k)\right]\cap A_n<
A_n.\] Renaming the elements in $\Omega$, we can assume $\gamma_1=\{1,\ldots,a\},\
\gamma_2\cup \ldots\cup\gamma_k=\{a+1,\ldots,n\}$ and, passing from $K$ to a proper
subgroup of $A_n$ containing $K$, we can assume \[K=\left[Sym\,\{1,\ldots,a\}\times Sym\,\
\{a+1,\ldots,n\}\right]\cap A_n\] with orbits $\{1,\ldots,a\}$ and $\{a+1,\ldots,n\}$ of
lengths $a$ and $b$ such that \[1\leq a\leq b\leq n-1,\ a+b=n.\]
Note that $a\leq \left[\frac{n}{2}\right];$ moreover no even permutation of the type
$[l_1;\ldots;l_r]$ with $l_i\geq n-a+1$ for some $i\in \{1,\ldots,r\}$ belongs to $K$,
because $b=n-a$ is the maximal length of an orbit for $K.$ Therefore $H$ contains at
least one permutation for each of these types. Observe that if we can identify a permutation
$\mu \in H$ of the type $[p\,;l_2;\ldots ;l_r]$ with $p$ prime and $p\nmid l_2,\ldots ,l_r,$ then
$\mu^{l_2\cdots l_r}\in H$ is a $p$-cycle. From now on our duty consists essentially in the
chase of these kinds of permutations and to do this we need to divide our argument into six
cases according to whether \[a\geq 6\,\hbox{ or }\, a=5,4,3,2,1.\]

\noindent {\em Case 1: $a\geq 6$}   

\noindent
Due to \ref{3.6}, it is sufficient to show that $H$ is primitive and contains a
$3$-cycle.

Let first $n$ be odd. Then $H$ contains permutations of the types: \[[n-2],\ [3;n-4],
\ [2;3;n-5].\] Because $H$ is transitive, $H_1$ contains a $(n-2)$-cycle and a
permutation of the type $[3;n-4]$ which has no fixed point on $\Omega-\{1\}.$ This
implies that $H_1$ is transitive on $\Omega-\{1\}$ and therefore $H$ is $2$-transitive
on $\Omega$ and, in particular, primitive.\\
If $3\nmid n-4,$ then a suitable power of a permutation in $H$ of type $[3;n-4]$ is a $3$-cycle.\\
If $3\mid n-4,$ then $3\nmid n-5$ and we get a $3$-cycle in $H$ as a power of a permutation
of the type $[2;3;n-5].$

Let now $n$ be even. Then $H$ contains permutations of the types:
\[[n-1],\ [3;n-3],\ [3;n-5].\]
$H$ is clearly $2$-transitive and then, as before, primitive.\\
If $3\nmid n,$ then a suitable power of a permutation of the type $[3;n-3]$ gives a $3$-cycle
in $H$.\\
If $3 \mid n,$ then $3 \nmid n-5$ and, we get $3$-cycle considering a power of a permutation
of the type $[3;n-5].$

~\\
\noindent {\em Case 2: $a=5$}   

\noindent
It is $b=n-5$ and $5\leq \left[\frac{n}{2}\right]$ leads to $n\geq 10.$ For $n\neq 12,$ it
is sufficient to show that $H$ is primitive and it contains a $3$-cycle or a $7$-cycle; the
case $n=12$ needs a particular argument.

Let $n$ be odd. Then $H$ contains permutations of the types: \[[n-2],\ [3;n-4].\]
Therefore $H_1$, containing a $(n-2)$-cycle and a fixed point free permutation, is
transitive. This means that $H$ is $2$-transitive and in particular primitive.\\
If $3\nmid n-1,$ we get a $3$-cycle in $H$ as a power of a permutation of the type $[3;n-4].$\\
If $3\mid n-1,$ then $n=3k+1$ with $k\geq 4$ even. We examine first the case $k>4,$ that
is $n>13:$ then the integer $\frac{n-3}{2}$ is greater than $5$ and  no
permutation of the type $[3;\frac{n-3}{2};\frac{n-3}{2}]$ lies in $K.$ On the other
hand $3\nmid n-3$ hence, if $\mu\in H$ is of that type, we get that
$\mu^{\frac{n-3}{2}} \in H$ is a $3$-cycle.\\ Next let $n=13:$
observe that, because $b=8$, $K$ contains no permutations of type $[2;4;7]$ and then
there exists $\mu \in H$ of the type $[2;4;7].$ Thus $\mu^4\in H$ is a $7$-cycle.

Let $n$ be even. Then $H$ contains permutations of the types: \[[n-1],\ [3;n-3].\]
This forces $H$ to be $2$-transitive and hence primitive.\\
If $3\nmid n,$ then $H$ contains a $3$-cycle. \\
If $3\mid n,$ then $n=6k$ with $k\geq 2.$ We explore first the case $k>2,$ that is
$n\geq 18:$ then the integer $\frac{n-4}{2}$ is greater than $5$ and no permutation
of the type $[3;\frac{n-4}{2};\frac{n-4}{2}]$ is in $K;$ on the other hand $3\nmid
n-4$ and, as usual, $H$ contains a $3$-cycle.\\
Next let $n=12.$ Then
\[K=\left[Sym\,\{1,2,3,4,5\}\times Sym\,\{6,7,8,9,10,11,12\}\right]\cap A_{12}\]
and $H$ contains permutations of the types:
\[[11],\ [9],\ [8;2],\ [8;4],\ [6;2;2;2],\ [3;3;3;3].\]
Then $H_1$ contains a $11$-cycle, $H_{1\,2}$ a permutation of the type $[9]$ and one
of the type $[8;2]$ and $H_{1\,2\,3}$ a permutation of the type $[9].$ This implies that $H$ is
$4$-transitive. Moreover we observe that $H_{1\,2\,3\,4}$ contains a
permutation of the type $[3;3]$ and one of the type $[4;4]$ as a power respectively of an
element of the type $[6;2;2;2]$ and of an element of the type $[8;2].$ Thus
$12\mid |H_{1\,2\,3\,4}|$ and hence
$2^5\cdot 3^4\cdot 5\cdot 11\mid |H|.$ On the other hand, by the well known Bochert's
result on the limitation of the index of a primitive subgroup of the symmetric group
( see \cite{hup67}, II,\,4.6 ), we have $|S_{12}:H|\geq 6!$ and thus
$|H|\leq 2^6\cdot 3^3\cdot 5\cdot 7\cdot 11.$ This gives
$$|H|=2^5\cdot 3^4\cdot 5\cdot 11\cdot k$$
with $1\leq k\leq 4.$\\
Yet we cannot have $k=3,$ otherwise $|H|_3=|A_{12}|_3$ and $H$
would contain a $3$-Sylow of $A_{12},$ hence also a $3$-cycle and, by \ref{3.6}, it would
coincide with $A_{12}.$ It follows that
\[|H|\in \{ 2^5\cdot 3^4\cdot 5\cdot 11,\, 2^6\cdot 3^4\cdot 5\cdot 11,\,
2^7\cdot 3^4\cdot 5\cdot 11\}.\]
Renaming the elements in $\Omega,$ we can assume $\sigma=(1\,2\,3\,4\,5\,6\,7\,8\,9\,\,10\,\,11)\in H$ and
$\gen{\sigma}\in Syl_{11}(H).$ It is easily observed that $N_{A_{12}}\gen{\sigma}=\gen{\sigma}
\rtimes \gen{(2\,5\,6\,10\,4)(3\,9\,\,11\,\,8\,7)}.$ Therefore we can have either
$N_H\gen{\sigma}=\gen{\sigma},$ or $N_H\gen{\sigma}=N_{A_{12}}\gen{\sigma}.$ Yet both these
possibilities are incompatible with $|H:N_H\gen{\sigma}|\equiv 1\,(\,\hbox{mod}\ 11\,).$ Namely, if
$N_H\gen{\sigma}=\gen{\sigma},$ we get $|H:N_H\gen{\sigma}|\in \{2^5\cdot 3^4\cdot 5,\,
2^6\cdot 3^4\cdot 5,\, 2^7\cdot 3^4\cdot 5\}$ while if $N_H\gen{\sigma}=N_{A_{12}}\gen{\sigma}$
we get $|H:N_H\gen{\sigma}|\in \{2^5\cdot 3^4,\, 2^6\cdot 3^4,\, 2^7\cdot 3^4\}$ and calculation
shows that neither of these numbers is congruent $1\ (\,\hbox{mod}\ 11\,).$ Therefore there
is no $H$ as required.

~\\
\noindent {\em Case 3: $a=4$}   

\noindent
As usual we examine separately the case $n$ odd and the case $n$ even, remembering that
the maximal length of an orbit for $K$ is $b=n-4.$ For $n\neq 9,$ we will show that $H$
is primitive and it contains a $3$-cycle or a $5$-cycle; the case $n=9$ needs a peculiar
argument.

Let $n$ be odd. Then $H$ contains permutations of the types: \[[n-2],\ [2;n-3]\] and
the usual argument show that  $H$ is primitive.\\
Let $n>9,$ then $K$ contains no permutation of the type $[2;3;n-5].$\\
If $3\nmid n-2,$ then a suitable power of an element of $H$ of the type $[2;3;n-5]$ is a
$3$-cycle.\\
If $3\mid n-2,$ then $3 \nmid n-3,\ n-3$ is even and $n=2+3k$ where $k\geq 3$ is odd.
If $k>3,$ we have $n>11$ and hence $\frac{n-3}{2}>4;$ this implies that $K$ does not
contain permutations of the type $[3;\frac{n-3}{2};\frac{n-3}{2}]$ and then $H$ contains
a $3$-cycle.\\
For $k=3,$ that is $n=11,$ we can observe that \[K=\left[Sym\,\{1,2,3,4\}\times
Sym\,\{5,6,7,8,9,10,11\}\right]\cap A_{11}\] contains no permutations of the type
$[3;3;5],$ hence $H$ contains at least a permutation of this type and hence
a $5$-cycle.\\
Finally let $n=9.$ Then
\[K=\left[Sym\,\{1,2,3,4\}\times Sym\,\{5,6,7,8,9\}\right]\cap A_{9}\]
contains no permutations of the types $[9],\ [7],\ [2;6],$ hence $H$ contains at
least a permutation for each of these types and consequently also a permutation of the type
$[3;3].$ It is clear that $H$ is primitive and, since no power of a $9$-cycle is of the type
$[3;3],$ we argue that $|H|_3\geq 3^3.$ But we cannot have $|H|_3> 3^3$ otherwise
$|H|_3=|A_9|_3$ and $H$ would contain a $3$-cycle. Therefore if $P\in Syl_3(H),$ we have
$|P|=27$ and
$P$ contains permutations of the types $[9]$ and $[3;3].$\footnote{We emphasize that $A_9$
really admits a subgroup $H$ with $|H|_3=3^3$ containing permutations of the types $[9],\
[7],\ [2;6].$ In fact you can consider $H=SL(2,8)\rtimes\gen{\alpha}$ where $\alpha$ is
the automorphism of $SL(2,8)$ which extend to the matrices the action of the Frobenius
automorphism of $GF(8)$ of degree $3,\ x\mapsto x^2;$ the natural action of this group on the
$9$ points of the projective line $\mathcal{P}(1,8)$ represents faithfully $H$ as a
subgroup of $A_9$ with the required properties.}\\
Renaming the elements in
$\Omega,$ we may assume $\mu=(123456789)\in P$ and there must exist $\beta\in P$ of the
type $[3;3],$
such that $P=\gen{\mu}\rtimes \gen{\beta}.$ This gives $P\leq N_{A_9}\gen{\mu}.$ But
$N_{A_9}\gen{\mu}={\gen{\mu}}\rtimes \gen{\alpha}$ with $\alpha=(235986)(47)$ and
$\mu^{\alpha}=(135792468)=\mu^2.$ It follows that
\[P=\gen{\mu,\ \beta: \mu^9=\beta^3=1,\ \mu^{\beta}=\mu^4},\] where we have put
$\beta=\alpha^2=(258)(396).$ Therefore $P$ is nilpotent of class $2,$ $P'=\gen{\mu^3}$
and each element in $P$ has a unique representation as $\mu^i\beta^j$ with $0\leq i\leq
8$ and $0\leq j\leq 2.$\\
We have $(\mu^i\beta^j)^3=\mu^{3i}\,\beta^{3j}\,[\beta^j,
\mu^i]^3=\mu^{3i},$ hence the elements of order $9$ in $P$, that is the $9$-cycles in $P,$
are exactly the $\mu^i\beta^j$ with $i\in I=\{1,2,4,5,7,8\}$ and $j\in\{0,1,2\}.$\\
Now we show that these elements are $A_9$-conjugate.\\ First of all we have
\[\mu^{\alpha^k}=\mu^{2^k}\quad \hbox{for}\quad 0\leq k\leq 5\]
and, since mod $9,$ the powers $2^k$ describe the elements of $I,$ the $9$-cycles
$\mu^i$ with $i\in I$ are $A_9$-conjugates.\\ Next consider $\gamma=(258)\in A_9$ and observe
that $\mu^{\gamma}=(153486729)=\mu \beta.$ Then
\[\mu^{\gamma\alpha^k}=(\mu\beta)^{\alpha^k}=\mu^{\alpha^k}\beta\quad\hbox{for}\quad
0\leq k\leq 5,\]
runs over all the $9$-cycles of the type $\mu^i\beta.$\\
Similarly, since $\gamma\in C_{A_9}(\beta),$
\[\mu^{\gamma^2\alpha^k}=(\mu\beta)^{\gamma\alpha^k}=(\mu^{\gamma}\beta)^{\alpha^k}=
(\mu\beta^2)^{\alpha^k}=\mu^{\alpha^k}\beta^2\quad\hbox{for}\quad 0\leq k\leq 5,\]
describes all the $9$-cycles of the type $\mu^i\beta^2.$\\
Since the $3$-Sylow subgroups of $H$ are $H$-conjugates, we deduce that the $9$-cycles
of $H$ are $A_9$-conjugates. Now, by \ref{3.4}, the $9$-cycles of $A_9$ split into two
conjugacy classes $\gamma_1,\ \gamma_2$ and calling $C$ the set of the $9$-cycles in
$H,$ we can assume $C\subseteq \gamma_1.$ But then picking any $\sigma\in \gamma_2,$
we get that $\sigma\notin \bigcup_{g\in {A_{9}}}K^g$ because $K$ contains no $9$-cycle
and, on the other hand, $\sigma\notin \bigcup_{g\in {A_{9}}}H^g$ because the $9$-cycles
in $\bigcup_{g\in {A_{9}}}H^g$ are given by $\bigcup_{g\in {A_{9}}}C^g
\subseteq \gamma_1,$ against the definition of a covering.

Let $n$ be even. Then $H$ contains permutations of the types: \[[n-1],\ [3;n-3]\] and it is
primitive.\\
If $3\nmid n,$ then $H$ contains a $3$-cycle.\\
If $3\mid n,$ then $n=6k$ with $k\geq 2.$ If $k>2,$ that is $n\geq 18,$ then $K$ contains
no permutations of the type $[2;3;5;n-10]$ and moreover $3\nmid n-10;$
thus there exists a $3$-cycle in $H.$\\
If $n=12,$ then \[K=\left[Sym\,\{1,2,3,4\}\times Sym\,\{5,6,7,8,9,10,11,12\}\right]
\cap A_{12}.\]
This implies that $H$ contains a permutation of type $[5;7]$ and hence a $5$-cycle.

~\\
\noindent {\em Case 4: $a=3$}   

\noindent
For $n\neq 10,$ it is sufficient to show that $H$ is $2$-transitive and
that it contains a $5$-cycle or, provided that $n>10,$ a $7$-cycle; the case $n=10,$
needs a special argument.

Let $n$ be odd. We consider first the case $n\neq 9,11.$ Then, because $n-5>3,$ $H$ contains permutations of the types:
\[[n-2],\ [4;n-5]\]
and it is primitive.\\
Because $n-5$ is even and $\frac{n-5}{2}>3$ for $n\neq 9,11,$ then
$H$ contains also a permutation of the type $[5;\frac{n-5}{2};\frac{n-5}{2}].$\\
If $5\nmid n,$ then we find a $5$-cycle in $H.$ \\
If $5\mid n,$ we observe that $5\nmid \,\frac{n-5}{2}-1,\,\frac{n-5}{2}
+1;$ moreover $\frac{n-5}{2}-1$ and $\frac{n-5}{2}+1$ have the same parity and therefore
there exists an even permutation of the type $[5;\frac{n-5}{2}-1;\frac{n-5}{2}+1].$ Yet, since
$n\geq 15\,$ implies $\,\frac{n-5}{2}-1 >3,$ no permutation of this type belongs to
$K$ and $H$ contains a $5$-cycle.\\
If $n=9,$ we observe that $K$ does not contain permutations of the type $[2;2;5]$
and we deduce that $H$ contains a $5$-cycle.\\ If $n=11,$ we note that $K$ contains no
permutation of the type $[2;2;7],$ hence $H$ contains a $7$-cycle.

Let $n$ be even. Then $H$ contains a $(n-1)$-cycle and it is primitive
while $K$ contains no permutation  of the type $[5;n-5].$ Hence, provided that
$5\nmid n,$ we have a $5$-cycle in $H$. On the other hand, if $5\mid n$ and $n\neq 10,$ then
we observe that $n-6$ is even and
$\frac{n-6}{2}>3;$ hence $K$ contains no permutation of the type $[5;\frac{n-6}{2};
\frac{n-6}{2}].$ Because $5\nmid\ \frac{n-6}{2},$ it turns out that $H$ contains a $5$-cycle.\\
Finally let $n=10.$ Then
\[K=\left[Sym\,\{1,2,3\}\times Sym\,\{4,5,6,7,8,9,10\}\right]\cap A_{10},\]
$H$ contains permutations of the types:
\[[9],\ [4;6]\]
and consequently also permutations of the type $[3;3].$  Since no power of a $9$-cycle
is of the type $[3;3],$ we argue that $|H|_3\geq 3^3.$ On the other hand, $H$ is clearly
primitive hence it contains no $3$-cycle and therefore $|H|_3<|A_{10}|_3.$ Thus,
if  $P\in Syl_3H,$  then $|P|=3^3$ and $P$ contains permutations of the types $[9]$ and
$[3;3].$ Now, the same
argument used in the case $a=4,\ n=9,$ enables us to assume $P=\gen{\mu}\rtimes\gen{\beta},$
where $\mu=(123456789)$ and $\beta$ is of the type $[3;3].$ Then we have
$P\leq N_{A_{10}}\gen{\mu}=\gen{\mu}\rtimes\gen{\alpha}$ where $\alpha=(235986)(47)$ and,
setting $\beta=\alpha^2,$ we get
$$P=\gen{\mu,\ \beta: \mu^9=\beta^3=1,\ \mu^{\beta}=\mu^4}$$
At this point observing that, by \ref{3.4}, the $9$-cycles in $A_{10}$ split into two conjugacy
classes, we can repeat word by word the reasoning developed in the case
$a=4,\ n=9$ reaching, as there, a contradiction.

~\\
\noindent
{\em Case 5: $a=2$}   

\noindent
It is enough to show that $H$ is primitive and that it contains a $3$-cycle.

Let $n$ be odd. Then $H_1$ contains permutations of the type:
 \[[3;n-4],\ [4;n-5].\]
This implies that $H_1$ is transitive on $\Omega-\{1\},$ otherwise it would have two orbits
of lengths $\{3,n-4\}=\{4,n-5\}$ which is impossible because $n\neq8.$
Therefore $H$ is primitive.\\
If $3\nmid n-1,$ then $H$ contains a $3$-cycle.\\
If $3\mid n-1,$ then $3\nmid \frac{n-3}{2}>2$ and $K$ contains no permutation of the type
$[3;\frac{n-3}{2};\frac{n-3}{2}]$ and therefore $H$ contains a $3$-cycle.

Let $n$ be even. Then $H$ contains permutations of the types: \[[n-1],\ [3;n-3]\] and is
primitive.\\
If $3\nmid n$ we obtain a $3$-cycle in $H$ in the usual way.\\
If $3\mid n,$ then $n=6k\geq 12$ and $K$ contains no permutation of the type $[3;4;n-8].$
Since $3\nmid n-8,$ we obtain a $3$-cycle in $H.$

~\\
\noindent
{\em Case 6: $a=1$}   

\noindent
In this case \[K=Alt\{2,\ldots,n\}\] and $\bigcup_{g\in {A_n}}K^g$ consists exactly of the
permutations of $A_n$ with at least a fixed point. Therefore if  $[l_1;\ldots;l_r]$ where
$l_i\geq 2$ are integers  such that  $\sum_{i=1}^r l_i=n,$  is the type of an even
permutation, then $H$ contains at least one permutation of this type.

Let $n$ be odd. Observe that for $k$ even $H$ contains permutations $\sigma_k$ of the types:
\[[\underbrace{2;\ldots;2}_k;n-2k],\]
provided that $0\leq k\leq \frac{n-3}{2},$ and for $k$ odd $H$ contains permutations
$\sigma_k$ of the types:
\[[\underbrace{2;\ldots;2}_{k-2};4;n-2k],\]
provided that $3\leq k\leq \frac{n-3}{2}.$\\
Since $\sigma_k^4\in H$ is a $(n-2k)$-cycle, we argue that $H$ contains at least
a $d$-cycle for each $d\in D=\{x\in \bfn: x\hbox{ odd, }3\leq x\leq n,\ x\neq n-2\}.$ But, by
the Bertrand's postulate, there exists a prime $p$ with $\frac{n+1}{2}\leq p \leq n-4$
and obviously $p\in D.$ Then \ref{3.7} implies that $H$ is primitive and, by \ref{3.6}, we
obtain $H=A_n.$

Let $n$ be even. Then $H$ contains permutations $\sigma_k$ of the types:
\[[k;n-k]\]
for each $k\in \bfn$ with $2\leq k\leq \frac{n}{2}.$\\
By \ref{3.7} and \ref{3.6}, it is sufficient to show that $H$ contains a $p$-cycle for
some $p$ prime with $\frac{n}{2}<p\leq n-3.$\\
Now, by the Bertrand's postulate, there exists such a prime $p$ and we can write
$p=n-k_0$ for some $3\leq k_0\leq \frac{n-2}{2}.$ Then, since $k_0<p,$ we get that
$\sigma_{k_0}^{k_0}\in H$ is a $p$-cycle.
\end{proof}
\bigskip

\section{The $(\ast\ast)$-coverable symmetric groups}

The main result of the previous section ( \ref{3.9} ) and the opportunity to build coverings
by intersection ( \ref{3.3} ), enable us to
extend immediately our investigation to the $(\ast\ast)$-coverable symmetric groups of
degree greater than $8.$ To complete
our analysis we  need essentially to clear what happens for $S_7$ and $S_8:$ we do this
re-echoing the methods used in the last section.

\begin{lem}\label{S_7}
$S_7$ and $S_8$ are not $(\ast\ast)$-coverable groups.
\end{lem}
\begin{proof}
Let $\{H^g,K^g:g\in S_7\}$ be a covering of $S_7$ and, without loss of generality,
let $H$ contain a $7$-cycle. Then, being $7$ a prime, $H$ is primitive and $K$ can not be
primitive either, otherwise one of them contains transpositions and hence, by \ref{3.6}, coincides with
$S_7.$ This implies that $K$ is not transitive and therefore we can assume
\[K=Sym\,\{1,\ldots,a\}\times Sym\,\{a+1,\ldots,7\}\]
with $1\leq a\leq 3$ minimal length of an orbit of $K.$\\
If $a=1,3,$ then $K$ contains no permutations of the type $[2;5]$ and hence  there
exists $\mu\in H$ of this type; consequently the transposition $\mu^5$ belongs to $H$
and we can again appeal to \ref{3.6} to conclude that $H=S_7.$\\
If $a=2,$ we have that $K=Sym\,\{1,2\}\times Sym\,\{3,4,5,6,7\}$ contains no permutations of
the types $[6]$ and $[3;4].$ Therefore $H$ contains at least a permutation of the type
$[6]$ and a permutation $\mu$ of the type $[3;4].$ Then $H\neq A_7,$ $\mu^4\in H$ is
a $3$-cycle and \ref{3.6} leads to the contradiction $H=S_7.$

Next let $\{H^g,K^g:g\in S_8\}$ be a covering of $S_8$ with $H$ containing a $7$-cycle.\\
If $H$ is not transitive, then we may assume $H=S_7$ and $K$ must contain all the types of
permutations with no fixed points. In particular $K$ contains at least a permutation of the
types $[8],\ [5;3]$ and $[3;3;2].$ Since $K$ is transitive and $5\mid\,|K|,$ from \ref{3.7},
we argue that $K$ is primitive. But because it contains a transposition, by \ref{3.6}, we
obtain that $K=S_8.$\\
Hence $H$ is transitive and indeed $2$-transitive.\\We observe that $H$ contains neither
$3$-cycles nor $5$-cycles. Namely if we assume the contrary, then by \ref{3.6}, we get $H=A_8;$
yet this implies that $K$ contains all the types of odd permutations and in particular $K$
contains a $8$-cycle and permutations of the types $[3;4]$ and $[6].$ Hence
$K$ is $2$-transitive and, since it contains transpositions, we deduce $K=S_8.$\\
In particular
$H$ contains no permutation of the types $[3;5],\ [3;4]$ and $ [2;5].$ Hence $K$ contains at
least a permutation of these types. Next observe that $K$ cannot be transitive otherwise $K_1$
contains permutations of the types $[3;4],$ $[2;5]$ and $K$ is $2$-transitive. Hence, since at
least a transposition belongs to $K$ we would have $K=S_8.$\\
This means that $K$ admits two orbits of lengths $3$ and $5$ and we may assume that
$$K=Sym\,\{1,2,3\}\times Sym\,\{4,5,6,7,8\}.$$
Then $H$ contains a $6$-cycle and is $3$-transitive. But it is well known that a
$k$-transitive group of degree $\,n\,$ with $\,k>\frac{n}{3}\,$ contains $A_n$
( see \cite{hup67}, p. $154$ ).
Then $H\geq A_8$ and, since $H$ contains also odd permutations, we argue that $H=S_8,$ a
contradiction.
\end{proof}
\begin{te}\label{sim}
$S_n$ is  $(\ast\ast)$-coverable if and only if $\,3\leq n\leq 6.$
\end{te}
\begin{proof}
Due to \ref{2.9}, \ref{3.4} and \ref{S_7}, what remains to show is that if $n\geq 9,$
then $S_n$ is not $(\ast\ast)$-coverable. Assume the contrary and let
$\delta=\{H^g,K^g:g\in S_n\}$ be a covering of $S_n$ where $n\geq 9.$ If $H$ and $K$
contain both some odd permutation, then we have $S_n=A_nK=A_nH$ and, by \ref{3.3},
$\delta$ defines by intersection a covering of $A_n$, against \ref{3.9}. Hence exactly
one among $H$ and $K$ is included in $A_n$: we can assume that $H=A_n,\ K\not\leq A_n.$
It follows that $K$ contains all the types of odd permutations and, in particular,
transpositions. Evidently, to reach a contra\-diction, it is sufficient to show that $K$ is
$2$-transitive. If $n$ is even, then $K$ contains a $n$-cycle and is transitive;
moreover $K_1$ containing a permutation of the type $[3;n-4]$ and a permutation of the
type $[n-2]$ is certainly transitive. If $n$ is odd, then $K$ contains a permutation of
the type $[n-1]$ and a permutation of the type $[3;n-3]$ with no fixed point; therefore
$K$ is transitive and, since $K_1$ contains a $(n-1)$-cycle, we get that $K$ is
$2$-transitive.
\end{proof}

\vspace*{3cm}
\noindent
\begin{minipage}[t]{4.4cm}
Address of the author:\\
~\\
DIMADEFAS \\
Via C. Lombroso 6/17 \\
I-50134 Firenze\\
Italy\\
E-mail:\\dbubbo@dmd.unifi.it
\end{minipage}
\hfill

\end{document}